\documentclass{article}%
\usepackage{amsmath}
\usepackage{graphicx}%
\usepackage{amsfonts}%
\usepackage{amssymb}
\newtheorem{theorem}{Theorem}

\newtheorem{corollary}[theorem]{Corollary}

\newtheorem{example}[theorem]{Example}

\newtheorem{lemma}[theorem]{Lemma}

\newtheorem{remark}[theorem]{Remark}

\newenvironment{proof}[1][Proof]{\textbf{#1.} }{\ \rule{0.5em}{0.5em}}

\begin{document}

\title{Lipschitz Flow-box Theorem}
\author{Craig Calcaterra and Axel Boldt\\{\small Department of Mathematics, Metropolitan State University, Saint Paul,
MN 55106}\\{\small craig.calcaterra@metrostate.edu, phone: (651) 793-1423, fax: (651) 793-1446}}
\maketitle

\begin{abstract}
A generalization of the Flow-box Theorem is proven. The assumption of a
$C^{1}$ vector field $f$ is relaxed to the condition that $f$ be locally
Lipschitz. The theorem holds in any Banach space.


\end{abstract}

\section{ Introduction}

The Flow-box Theorem for smooth vector fields states that the dynamic near a
non-equilibrium point is qualitatively trivial, i.e., topologically conjugate
with translation. Near a nondegenerate equilibrium point, linearizing the
vector field by differentiation allows a relatively simple characterization of
almost all possible local dynamics. These two results characterize the local
behavior of solutions for smooth nondegenerate vector fields. A natural
follow-up question is, ``What dynamics are possible under nonsmooth conditions?''

To be more specific, the Flow-box Theorem (also called the ``Straightening-out
Theorem'' or the ``Local Linearization Lemma'') applies to autonomous,
first-order differential equations, i.e.,%
\begin{equation}
x^{\prime}\left(  t\right)  =f\left(  x\right(  t))\text{.} \label{diffEq}%
\end{equation}
$f$ typically is a vector field on a manifold. For local questions such as
ours, it is enough to study the case of a map $f:X\rightarrow X$ where
$X=\mathbb{R}^{n}$ or some other Banach space. A \textbf{solution} to $\left(
\ref{diffEq}\right)  $ with initial condition $x_{0}\in X$ is a curve
$x:I\rightarrow X$ where $I$ is an open subinterval of $\mathbb{R}$ containing
$0$, $x\left(  0\right)  =x_{0},$ and which satisfies $\left(  \ref{diffEq}%
\right)  $ for all $t\in I$.

The traditional Flow-box Theorem asserts that if $f$ is a $C^{1}$ vector field
and $x_{0}\in X$ is not an equilibrium, i.e., $f(x_{0})\neq0$, then there is a
diffeomorphism which transfers the vector field near $x_{0}$ to a constant
vector field. In other words, the local flow of $f$ is conjugate via
diffeomorphism to translation.

The Picard-Lindel\"{o}f Theorem\footnote{Also known as The Cauchy-Lipschitz
Theorem, The Fundamental Theorem of Differential Equations, or the Local
Existence and Uniqueness Theorem. It is proven, e.g., in \cite[p.
188]{Abraham}.}, stated below, guarantees a unique solution $\sigma_{x}$
exists for every initial condition $x\in X$ if $f$ is locally
Lipschitz-continuous. The continuous dependence of solutions on initial
conditions (Lemma \ref{LipEst}, below) is also assured when $f$ is Lipschitz
continuous. For these nonsmooth vector fields, is the dynamic near
non-equilibria still qualitatively trivial? I.e., does the Flow-box Theorem
still hold when we drop the $C^{1}$ condition on $f$? Yes and no.

A transferring diffeomorphism need not exist if $f$ is merely Lipschitz
(Example \ref{diffEx} below). So the Flow-box Theorem does not trivially
extend to the Lipschitz case. The natural next hypothesis is that the Flow-box
Theorem for Lipschitz vector fields might work if we use a transferring
\textit{lipeomorphism }(a bijective Lipschitz map whose inverse is also
Lipschitz). Lipeomorphisms cannot transfer vector fields, but they can still
provide a conjugacy between two dynamics. The result of this paper, Theorem
\ref{VFLinear}, is that for every non-equilibrium of a Lipschitz vector field
there exists a local conjugacy to a constant vector field via lipeomorphism.
Therefore the topological conjugacy with translation still holds when a vector
field is not differentiable.

Roughly, the trick in constructing the flow box is to track solutions to a
hyperplane transverse to the vector $f\left(  x_{0}\right)  $. The traditional
proofs then employ the Implicit Function Theorem or Inverse Function Theorem
requiring differentiability. For merely Lipschitz conditions, generalizations
of those theorems exist, but do not help when checking the transferring map is
Lipschitz. We rely on the Picard-Lindel\"{o}f Theorem and Lipschitz continuous
dependence on initial conditions to finish the proof. We do not make use of
Rademacher's Theorem which says a Lipschitz map is almost everywhere differentiable.

For manifolds the Flow-box Theorem states that for any $C^{1}$ vector field
with $f(x)\neq0$ there is a chart around $x$ on which $f$ is constant. Proofs
for $C^{\infty}$ Banach manifolds can be found in \cite{Abraham} or
\cite{Lang}. The results of this paper are easily ported to this context: a
vector field is called locally Lipschitz continuous if it is locally Lipschitz
in one chart (and therefore all charts).

Thus the local qualitative characterization of dynamical systems under
Lipschitz conditions reduces to the study of equilibria. This question has
already been broached, as dynamics with nonsmooth vector fields has enjoyed
some popularity in the last few decades in control theory. Discontinuous
vector fields have been analyzed with a host of different approaches: see for
instance \cite{AubinFrank}, \cite{Broucke}, \cite{Clarke}, \cite{Deimling},
\cite{Kunze}. Even for the less extreme case of Lipschitz continuous vector
fields, the analysis of equilibria is ever more complicated than the smooth
non-degenerate hyperbolic case.

Interesting related results have been obtained in \cite{Abergel} concerning
Lyapunov exponents for systems generated by Lipschitz vector fields, and in
\cite{RampSuss} where the Lie bracket is generalized to Lipschitz vector
fields. The Flow-box Theorem is the base case for Frobenius' Theorem on the
equivalence of involutive and integrable distributions. \cite{Simic} presents
a generalization of Frobenius' Theorem for Lipschitz vector fields.

We finish with some examples that explore how the quality of continuity of a
vector field ($C^{0}$ vs. Lipschitz vs. $C^{1}$) are related to the quality of
continuity of the transferring map in the flow box setting.

\section{Lipschitz Flow-box Theorem}

A \textbf{Banach space} is a normed linear space, complete in its norm.

A map $f:U\rightarrow V$ between subsets of Banach spaces is
\textbf{Lipschitz} if there exists $K>0$ such that%
\[
\left\|  f\left(  x_{1}\right)  -f\left(  x_{2}\right)  \right\|  \leq
K\left\|  x_{1}-x_{2}\right\|
\]
for all $x_{1},x_{2}\in U$. A \textbf{lipeomorphism} is an invertible
Lipschitz map whose inverse is also Lipschitz (i.e., slightly stronger than a
homeomorphism). A \textbf{vector field} on a Banach space $X$ is a map
$f:U\rightarrow X$ where $U\subset X$. A \textbf{solution} to a vector field
$f$ with \textbf{initial condition} $x$ is a curve $\sigma_{x}:I\rightarrow U$
defined on an open interval $I$ containing $0$ such that $\sigma_{x}\left(
0\right)  =x$ and $\sigma_{x}^{\prime}\left(  t\right)  =f\left(  \sigma
_{x}\left(  t\right)  \right)  $ for all $t\in I$.

Denote the open ball in $X$ about $x_{0}\in X$ with radius $r$ by%
\[
B\left(  x_{0},r\right)  :=\left\{  x\in X:\left\|  x-x_{0}\right\|
<r\right\}  .
\]

\begin{theorem}
[Picard-Lindel\"{o}f]\label{Picard-L}Let $X$ be a Banach space, $x_{0}\in X$,
and $r>0$. Assume $f:B\left(  x_{0},r\right)  \rightarrow X$ is a Lipschitz
vector field, and let $M$ be such that $\left\|  f\left(  x\right)  \right\|
\leq M$ for all $x\in B\left(  x_{0},r\right)  $. Then there exists a unique
solution to $f$ with initial condition $x_{0}$ defined on $\left(
-\frac{r}{M},\frac{r}{M}\right)  $.
\end{theorem}

\begin{proof}
See, e.g., \cite[p. 188]{Abraham} for the idea.
\end{proof}

The following well-known result (also given in \cite[p. 189]{Abraham}) is used
in the proof of the main theorem.

\begin{lemma}
[Continuous dependence on initial conditions]\label{LipEst}Let $f$ be a
Lipschitz vector field with constant $K$ defined on an open subset of a Banach
space. Let $\sigma_{x}$ and $\sigma_{y}$ be solutions to $f$ for initial
conditions $x$ and $y$ with interval $I\ni0$ contained in their common
domains. Then%
\[
\left\|  \sigma_{x}\left(  t\right)  -\sigma_{y}\left(  t\right)  \right\|
\leq\left\|  x-y\right\|  e^{K\left|  t\right|  }%
\]
for all $t\in I$.
\end{lemma}

As a consequence of the two previous results, for any Lipschitz vector field
$f:U\rightarrow X$ on an open set $U$, near any point $x\in U$ there exists a
local flow. Specifically:

\begin{corollary}
Let $f:U\rightarrow X$ be a Lipschitz vector field on a Banach space $X$. Let
$x\in U$. There exists a neighborhood $W$ of $x$ in $U$, a number $\epsilon
>0$, and a map $F:W\times\left(  -\epsilon,\epsilon\right)  \rightarrow X$
such that

1. $\frac{d}{dt}F\left(  x,t\right)  =f\left(  F\left(  x,t\right)  \right)  $

2. $F\left(  x,0\right)  =x$

3. $F\left(  x,s+t\right)  =F\left(  F\left(  x,s\right)  ,t\right)  $.
\end{corollary}

$F$ is called a \textbf{local flow} of $f$ near $x.$

Two vector fields $f_{1}:U_{1}\rightarrow X$ and $f_{2}:U_{2}\rightarrow X$
are called \textbf{locally topologically conjugate} near $x_{1}\in U_{1}$ and
$x_{2}\in U_{2}$ if there exist open neighborhoods $W_{1}$ of $x_{1}$ and
$W_{2}$ of $x_{2}$ and a homeomorphism $\phi:W_{1}\rightarrow W_{2}$ with
$\phi\left(  x_{1}\right)  =x_{2}$ such that a curve $\sigma:I\rightarrow
W_{1}$ is a solution to $f_{1}$ if and only if $\phi\circ\sigma:I\rightarrow
W_{2}$ is a solution to $f_{2}$. Less formally we require the local flows
$F_{1}$ and $F_{2}$ to satisfy%
\[
\phi\left(  F_{1}\left(  x,t\right)  \right)  =F_{2}\left(  \phi\left(
x\right)  ,t\right)
\]
for all $x\in W_{1}$ and $\left|  t\right|  $ sufficiently small.

Now we are ready to state and prove the main result of the paper.

\begin{theorem}
[Flow box]\label{VFLinear}Let $X$ be a Banach space with open subset $U$. Let
$f:U\rightarrow X$ be a Lipschitz vector field. Let $z\in X$ be nonzero and
let $g:X\rightarrow X$ be the constant vector field $g\left(  x\right)  =z$.
Then for any point $x_{1}\in U$ with $f\left(  x_{1}\right)  \neq0$, $f$ and
$g$ are locally topologically conjugate near $x_{1}$ and $x_{2}:=0$.

The homeomorphism which gives the conjugacy is a lipeomorphism.
\end{theorem}

\begin{description}
\item[Remark:] In the usual formulation of the flow-box theorem in
$\mathbb{R}^{n}$, the constant vector field $g$ is usually chosen as $g\left(
x\right)  =\left(  1,0,...,0\right)  $.
\end{description}

\begin{proof}
(Outline) We may assume $x_{1}=0$ and $f\left(  0\right)  =z$. First specify a
hyperplane $\Pi$ transverse to the vector $f\left(  x_{1}\right)  $. Then for
$x\in X$ near $x_{1}$ track solutions $\sigma_{x}$ back to the plane $\Pi$.
Define $t_{x}$ to be the value of $t$ such that $\sigma_{x}\left(  -t\right)
\in\Pi$ and define $p_{x}:=\sigma_{x}\left(  -t_{x}\right)  $. Then the
transferring lipeomorphism $\phi$ we seek is $\phi\left(  x\right)
:=p_{x}+t_{x}z$.
\begin{figure}
[ptbh]
\begin{center}
\includegraphics[
height=1.9026in,
width=4.708in
]%
{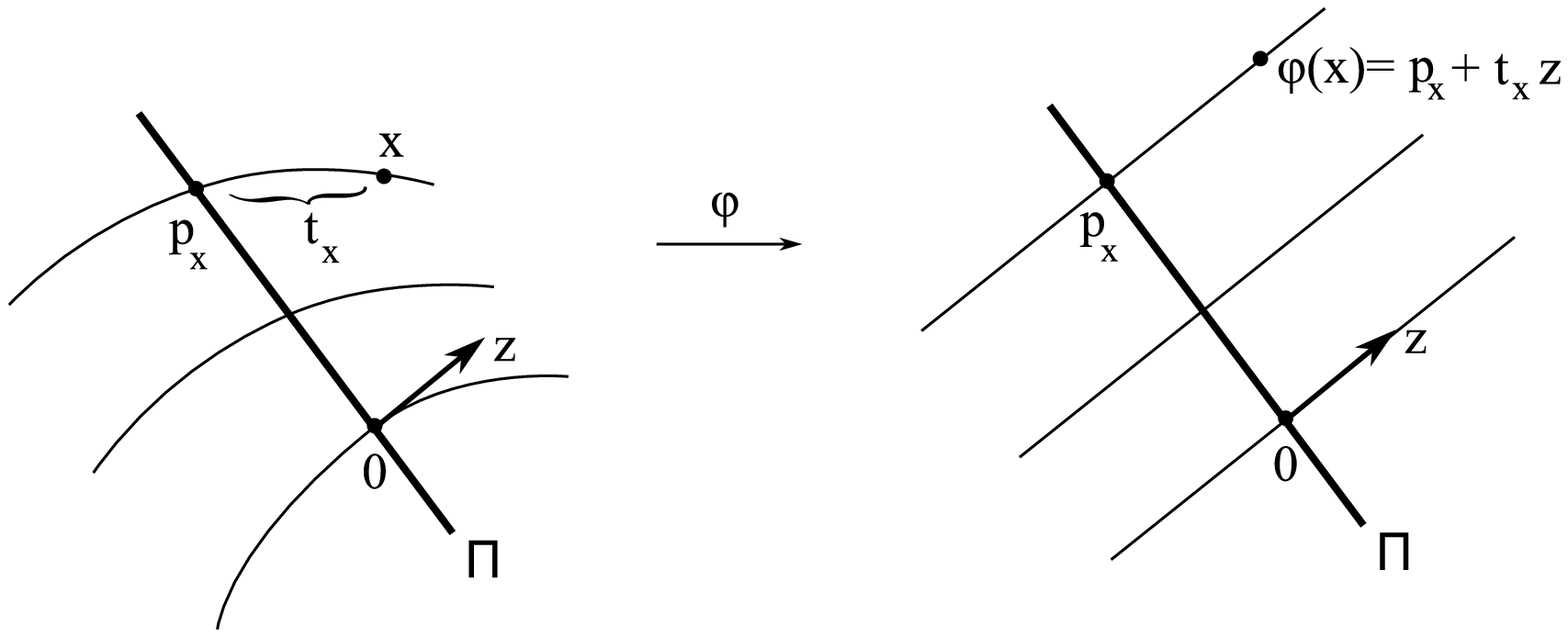}%
\end{center}
\end{figure}
We check that $\phi$ supplies the conjugacy and finally that $\phi$ is a
lipeomorphism, but first we must find a neighborhood $W_{1}$ of $x_{1}$ on
which $\phi$ is well-defined. Finding $W_{1}$ is the difficult part of the
proof and this culminates at $\left(  \ref{px}\right)  $ below.

(Details) First we verify we may assume without loss of generality that
$x_{1}=0$ and $f\left(  0\right)  =z$ and $\left\|  z\right\|  =1$. The
translation to $0$ and the dilation to norm one are obvious. If $z$ and
$f\left(  0\right)  =:y$ are linearly independent, then the intermediate
transferring diffeomorphism $A$ is a little more difficult to construct.
Consider the function $\psi$ from the subspace spanned by $y$ and $z$ to
$\mathbb{R}$ given by $\psi\left(  ay+bz\right)  :=b-a$. Extend $\psi$ to a
continuous linear functional $\overline{\psi}$on $X$ with the Hahn-Banach
Theorem. Then $A:X\rightarrow X$ given by $A\left(  x\right)  :=x+\overline
{\psi}\left(  x\right)  \left(  y-z\right)  $ is its own inverse and does the job.

We make several successive refinements of a neighborhood of $0$ in
constructing $W_{1}$ and the lipeomorphism $\phi:W_{1}\rightarrow W_{2}$.

By the Hahn-Banach Theorem there exists a continuous $\mathbb{R}$-linear map
$\chi:X\rightarrow\mathbb{R}$ with $\chi\left(  z\right)  =1$ and $\left|
\chi\left(  x\right)  \right|  \leq\left\|  x\right\|  $ for all $x\in X$. By
the continuity of $f$ there is some $r_{1}>0$ such that
\[
\chi\left(  f\left(  x\right)  \right)  >\frac{1}{2}\text{\hspace
{0.2in}and\hspace{0.2in}}\left\|  f\left(  x\right)  \right\|  <2
\]
for all $x$ in the open ball $B\left(  0,r_{1}\right)  $. We can also assume
$r_{1}$ is chosen so the Lipschitz condition for $f$ is met on all of
$B\left(  0,r_{1}\right)  $ with constant $K$.

By Theorem \ref{Picard-L} and the triangle inequality, for each $x\in B\left(
0,r_{1}/2\right)  $ a unique solution to $f$ with initial condition $x$
exists, $\sigma_{x}:\left(  -T,T\right)  \rightarrow B\left(  0,r_{1}\right)
$ where $T:=r_{1}/4$.

Denote the hyperplane in $X$ which is the kernel of $\chi$ by%
\[
\Pi:=\left\{  x:\chi\left(  x\right)  =0\right\}  \text{.}%
\]
Define
\[
R:=\underset{p\in B\left(  0,r_{1}/2\right)  \cap\Pi}{\cup}\sigma_{p}\left(
\left(  -T,T\right)  \right)  .
\]
With $r_{2}=\min\left\{  \frac{r_{1}}{10},\frac{T}{2}\right\}  $ define
\[
W_{1}:=B\left(  0,r_{2}\right)  \text{.}%
\]

Next we show $W_{1}\subset R$. This will then guarantee that every point $x\in
W_{1}$ is obtained by following a solution $\sigma_{p}$ with initial condition
$p\in\Pi$. For $x\in W_{1}$ we know $\sigma_{x}\left(  \left(  -T,T\right)
\right)  \subset B\left(  0,r_{1}\right)  $. Then
\[
\left(  \chi\circ\sigma_{x}\right)  ^{\prime}\left(  t\right)  =\chi\left(
\sigma_{x}^{\prime}\left(  t\right)  \right)  =\chi\left(  f\left(  \sigma
_{x}\left(  t\right)  \right)  \right)  >\frac{1}{2}%
\]
for $-T<t<T$. Further
\[
\left|  \chi\circ\sigma_{x}\left(  0\right)  \right|  =\left|  \chi\left(
x\right)  \right|  \leq\left\|  x\right\|  <r_{2}\text{.}%
\]
Thus there exists a unique $t\in\left(  -2r_{2},2r_{2}\right)  $ such that
$\chi\left(  \sigma_{x}\left(  t\right)  \right)  =0$, i.e., $\sigma
_{x}\left(  t\right)  \in B\left(  0,r_{1}\right)  \cap\Pi$. Furthermore the
speed of $\sigma_{x}$ is less than $2$ so that the distance from $x$ to
$\sigma_{x}\left(  t\right)  $ has $4r_{2}$ as an upper bound. This follows
since $\left\|  \sigma_{x}\left(  0\right)  -\sigma_{x}\left(  t\right)
\right\|  \leq\left|  \int_{0}^{t}\left\|  \sigma_{x}^{\prime}\left(
s\right)  \right\|  ds\right|  $. Thus the distance from $0$ to $\sigma
_{x}\left(  t\right)  $ is less than $5r_{2}\leq\frac{r_{1}}{2}$ and so
$\sigma_{x}\left(  t\right)  \in B\left(  0,r_{1}/2\right)  \cap\Pi$. Due to
the uniqueness of solutions, $\sigma_{\sigma_{x}\left(  t\right)  }\left(
-t\right)  =x$ so that $x\in R$ and the claim is proven.

Thanks to the previous paragraph we can now define $\phi$. For each $x\in
W_{1}$ since $\chi\circ\sigma_{x}$ is strictly increasing on $\left(
-T,T\right)  $, there exists a unique $t_{x}\in\left(  -T,T\right)  $ such
that
\begin{equation}
\sigma_{x}\left(  -t_{x}\right)  \in B\left(  0,r_{1}/2\right)  \cap
\Pi\text{.}\label{px}%
\end{equation}
Define
\[
p_{x}:=\sigma_{x}\left(  -t_{x}\right)
\]
and define $\phi$ by
\[
\phi\left(  x\right)  :=p_{x}+t_{x}z\in X
\]
for each $x\in W_{1}$.

Let $F$ be the local flow of $f$ near $x_{1}=0$ and let $G$ be the flow of
$g$. I.e., $F\left(  x,t\right)  :=\sigma_{x}\left(  t\right)  $ and $G\left(
x,t\right)  =x+tz$. To demonstrate the conjugacy we show $\phi\left(  F\left(
x,t\right)  \right)  =G\left(  \phi\left(  x\right)  ,t\right)  $ whenever
$x\in W_{1}$ and $\left|  t\right|  $ is sufficiently small. Notice that the
definitions of $p_{x}$ and $t_{x}$ above give%
\[
p_{F\left(  x,t\right)  }=p_{x}\text{\hspace{0.3in}and\hspace{0.3in}%
}t_{F\left(  x,t\right)  }=t_{x}+t
\]
so that
\[
\phi\left(  F\left(  x,t\right)  \right)  =p_{F\left(  x,t\right)
}+t_{F\left(  x,t\right)  }z=p_{x}+\left(  t_{x}+t\right)  z=G\left(
\phi\left(  x\right)  ,t\right)  \text{.}%
\]

$\phi$ is 1-1. To see this suppose $\phi\left(  x\right)  =\phi\left(
y\right)  $. Then $p_{x}-p_{y}=\left(  t_{y}-t_{x}\right)  z$. Applying $\chi$
yields $t_{x}=t_{y}$ so that $p_{x}=p_{y}$. By the uniqueness of solutions to
$f$ we get $x=\sigma_{p_{x}}\left(  t_{x}\right)  =\sigma_{p_{y}}\left(
t_{y}\right)  =y$.

To show Lipschitz continuity we will use Lemma \ref{LipEst}. Pick $x,y\in
W_{1}$. Since $\left(  \chi\circ\sigma_{x}\right)  ^{\prime}\left(  t\right)
>\frac{1}{2}$ for all $t\in\left(  -T,T\right)  $,%
\begin{align*}
&  \left|  t_{x}-t_{y}\right| \\
&  \leq2\left|  \left(  \chi\circ\sigma_{x}\right)  \left(  -t_{x}\right)
-\left(  \chi\circ\sigma_{x}\right)  \left(  -t_{y}\right)  \right| \\
&  \leq2\left(  \left|  \left(  \chi\circ\sigma_{x}\right)  \left(
-t_{x}\right)  -\left(  \chi\circ\sigma_{y}\right)  \left(  -t_{y}\right)
\right|  +\left|  \left(  \chi\circ\sigma_{y}\right)  \left(  -t_{y}\right)
-\left(  \chi\circ\sigma_{x}\right)  \left(  -t_{y}\right)  \right|  \right)
\\
&  \leq2\left(  \left|  \chi\left(  p_{x}\right)  -\chi\left(  p_{y}\right)
\right|  +\left\|  \sigma_{y}\left(  -t_{y}\right)  -\sigma_{x}\left(
-t_{y}\right)  \right\|  \right) \\
&  \leq2\left(  0+\left\|  x-y\right\|  e^{K\left|  t_{y}\right|  }\right)
\end{align*}
Next, using the bound on speed $\left\|  f\left(  x\right)  \right\|  <2$
gives%
\begin{gather*}
\left\|  p_{x}-p_{y}\right\|  =\left\|  \sigma_{x}\left(  -t_{x}\right)
-\sigma_{y}\left(  -t_{y}\right)  \right\| \\
\leq\left\|  \sigma_{x}\left(  -t_{x}\right)  -\sigma_{x}\left(
-t_{y}\right)  \right\|  +\left\|  \sigma_{x}\left(  -t_{y}\right)
-\sigma_{y}\left(  -t_{y}\right)  \right\| \\
\leq2\left|  t_{x}-t_{y}\right|  +\left\|  x-y\right\|  e^{K\left|
t_{y}\right|  }\leq\left\|  x-y\right\|  5e^{K\left|  t_{y}\right|  }\text{.}%
\end{gather*}
Since $\left|  t_{y}\right|  <T,$ defining $K_{\phi}:=7e^{KT}$ gives%
\begin{align*}
\left\|  \phi\left(  x\right)  -\phi\left(  y\right)  \right\|   &  =\left\|
p_{x}+t_{x}z-\left(  p_{y}+t_{y}z\right)  \right\|  =\left\|  \left(
t_{x}-t_{y}\right)  z+\left(  p_{x}-p_{y}\right)  \right\| \\
&  \leq\left|  t_{x}-t_{y}\right|  +\left\|  p_{x}-p_{y}\right\|  \leq
K_{\phi}\left\|  x-y\right\|  \text{.}%
\end{align*}

Now we show $\phi^{-1}$ is Lipschitz. Pick $u=p_{x}+t_{x}z=\phi\left(
x\right)  $ and $v=p_{y}+t_{y}z=\phi\left(  y\right)  $ then%
\begin{align*}
&  \left\|  \phi^{-1}\left(  u\right)  -\phi^{-1}\left(  v\right)  \right\|
=\left\|  x-y\right\|  =\left\|  \sigma_{p_{x}}\left(  t_{x}\right)
-\sigma_{p_{y}}\left(  t_{y}\right)  \right\| \\
&  \leq\left\|  \sigma_{p_{x}}\left(  t_{x}\right)  -\sigma_{p_{y}}\left(
t_{x}\right)  \right\|  +\left\|  \sigma_{p_{y}}\left(  t_{x}\right)
-\sigma_{p_{y}}\left(  t_{y}\right)  \right\|  \text{.}%
\end{align*}
Using Lemma \ref{LipEst} again, we get%
\[
\left\|  \sigma_{p_{x}}\left(  t_{x}\right)  -\sigma_{p_{y}}\left(
t_{x}\right)  \right\|  \leq\left\|  p_{x}-p_{y}\right\|  e^{K\left|
t_{x}\right|  }%
\]
and the bound on speed $\left\|  f\left(  x\right)  \right\|  <2$ gives%
\[
\left\|  \sigma_{p_{y}}\left(  t_{x}\right)  -\sigma_{p_{y}}\left(
t_{y}\right)  \right\|  \leq2\left|  t_{x}-t_{y}\right|  \text{.}%
\]
Define the projection $\pi:X\rightarrow\Pi$ along $z$ by%
\[
\pi\left(  q\right)  :=q-\chi\left(  q\right)  z\text{.}%
\]
This is a linear map and continuous since%
\[
\left\|  \pi\left(  q\right)  \right\|  \leq\left\|  q\right\|  +\left|
\chi\left(  q\right)  \right|  \left\|  z\right\|  \leq2\left\|  q\right\|
\text{.}%
\]
Then%
\[
\left\|  p_{x}-p_{y}\right\|  =\left\|  \pi\left(  u\right)  -\pi\left(
v\right)  \right\|  \leq2\left\|  u-v\right\|
\]
and
\[
\left|  t_{x}-t_{y}\right|  =\left|  \chi\left(  u\right)  -\chi\left(
v\right)  \right|  \leq\left\|  u-v\right\|  \text{.}%
\]
Define $K_{\phi^{-1}}:=2+2e^{KT}$. Combining the above estimates gives%
\begin{align*}
&  \left\|  \phi^{-1}\left(  u\right)  -\phi^{-1}\left(  v\right)  \right\| \\
&  \leq\left\|  p_{x}-p_{y}\right\|  e^{K\left|  t_{x}\right|  }+2\left|
t_{x}-t_{y}\right| \\
&  \leq K_{\phi^{-1}}\left\|  u-v\right\|
\end{align*}
so $\phi^{-1}$ is Lipschitz.

Finally since $X$ might be infinite dimensional, we must check $W_{2}%
:=\phi\left(  W_{1}\right)  $ is open. Let $p_{x}+t_{x}z=\phi\left(  x\right)
\in W_{2}$ for $x\in W_{1}$. Since $W_{1}$ is open there exists $s_{1}>0$ such
that $B\left(  x,s_{1}\right)  \subset W_{1}$. Since $t_{x}\in\left(
-T,T\right)  $, $s_{2}:=\min\left\{  T-\left|  t_{x}\right|  ,\frac{s_{1}}%
{4}\right\}  >0$. Then pick $s_{3}>0$ such that $B\left(  p_{x},s_{3}\right)
\subset B\left(  0,\frac{r_{1}}{2}\right)  $ and such that for all $p\in
B\left(  p_{x},s_{3}\right)  $ we have $\left\|  \sigma_{p}\left(
t_{x}\right)  -\sigma_{p_{x}}\left(  t_{x}\right)  \right\|  <\frac{s_{1}}{2}$
(using Lemma \ref{LipEst}). Then with $s_{4}:=\min\left\{  s_{2},\frac{s_{3}%
}{2}\right\}  >0$ we have $B\left(  \phi\left(  x\right)  ,s_{4}\right)
\subset W_{2}$. To see this notice any member of $B\left(  \phi\left(
x\right)  ,s_{4}\right)  $ may be written uniquely as $p+tz$ for some $p\in
\Pi$ and $t\in\mathbb{R}$. Then
\begin{align*}
\left|  t-t_{x}\right|   &  =\left|  \chi\left(  \left[  p_{x}+t_{x}z\right]
-\left[  p+tz\right]  \right)  \right| \\
&  \leq\left\|  \left[  p_{x}+t_{x}z\right]  -\left[  p+tz\right]  \right\|
<s_{4}\leq s_{2}%
\end{align*}
and%
\begin{align*}
\left\|  p-p_{x}\right\|   &  =\left\|  \pi\left(  \left[  p_{x}%
+t_{x}z\right]  -\left[  p+tz\right]  \right)  \right\| \\
&  \leq2\left\|  \left[  p_{x}+t_{x}z\right]  -\left[  p+tz\right]  \right\|
<2s_{4}\leq s_{3}\text{.}%
\end{align*}
Then%
\begin{align*}
\left\|  \sigma_{p}\left(  t\right)  -x\right\|   &  =\left\|  \sigma
_{p}\left(  t\right)  -\sigma_{p_{x}}\left(  t_{x}\right)  \right\| \\
&  \leq\left\|  \sigma_{p}\left(  t\right)  -\sigma_{p}\left(  t_{x}\right)
\right\|  +\left\|  \sigma_{p}\left(  t_{x}\right)  -\sigma_{p_{x}}\left(
t_{x}\right)  \right\|  <\tfrac{s_{1}}{2}+\tfrac{s_{1}}{2}=s_{1}%
\end{align*}
so $\sigma_{p}\left(  t\right)  \in W_{1}$ and therefore $\phi\left(
\sigma_{p}\left(  t\right)  \right)  =p+tz\in\phi\left(  W_{1}\right)  =W_{2}$.
\end{proof}

\begin{remark}
The Hahn-Banach theorem is essential for this proof. On a Hilbert space or
$\mathbb{R}^{n}$ with arbitrary norm, however, an obvious modification (in the
definition of $\Pi$) yields a proof which does not rely on the Axiom of Choice.
\end{remark}

\begin{example}
\label{diffEx}The lipeomorphism constructed in the proof is not necessarily
differentiable. If a vector field $f:\mathbb{R}^{n}\rightarrow\mathbb{R}^{n}$
is smooth then the traditional Flow-box Theorem states that near a
non-equilibrium $x_{0}$ there exists a diffeomorphism $\phi$ from a
neighborhood of $x_{0}$ to a neighborhood of the origin such that $\phi_{\ast
}f=g$, i.e., $d\phi_{x}\left(  f\left(  x\right)  \right)  =g\left(
\phi\left(  x\right)  \right)  $, where $g$ is the constant vector field
$g\left(  x\right)  :=\left(  1,0,...,0\right)  $. We show in this example
that if $f$ is not smooth, there may not exist such a diffeomorphism $\phi$.

Consider the Lipschitz vector field $f:\mathbb{R}^{2}\mathbb{\rightarrow
R}^{2}$ given by $f\left(  x,y\right)  =\left(  1+\left|  y\right|  ,0\right)
$. If a transferring diffeomorphism $\phi$ did exist near $x_{0}=\left(
0,0\right)  $, then with $\psi:=\phi^{-1}$ we would have%
\begin{align*}
\left(  d\psi\right)  _{\left(  x,y\right)  }\left(  g\left(  x,y\right)
\right)   &  =f\circ\psi\left(  x,y\right)  \\
\left[
\begin{array}
[c]{cc}%
\psi_{1x}\left(  x,y\right)   & \psi_{1y}\left(  x,y\right)  \\
\psi_{2x}\left(  x,y\right)   & \psi_{2y}\left(  x,y\right)
\end{array}
\right]
\genfrac{[}{]}{0pt}{0}{1}{0}%
&  =%
\genfrac{[}{]}{0pt}{0}{\psi_{1x}\left(  x,y\right)  }{\psi_{2x}\left(
x,y\right)  }%
=%
\genfrac{[}{]}{0pt}{0}{1+\left|  \psi_{2}\left(  x,y\right)  \right|  }{0}%
\end{align*}
therefore there are functions $u$ and $v$ such that%
\begin{align*}
\psi_{2}\left(  x,y\right)   &  =u\left(  y\right)  \\
\psi_{1x}\left(  x,y\right)   &  =1+\left|  u\left(  y\right)  \right|  \\
\psi_{1}\left(  x,y\right)   &  =x+x\left|  u\left(  y\right)  \right|
+v\left(  y\right)
\end{align*}
Since $\phi\left(  0,0\right)  =\left(  0,0\right)  $ we must have $u\left(
0\right)  =0$,  and since $\psi_{2x}\left(  x,y\right)  =0$ we cannot have
$u^{\prime}\left(  y\right)  =0$ lest the Jacobian be singular. Therefore
there is no neighborhood of the origin in which $\psi_{1}$ is differentiable.
\end{example}

\begin{remark}
In one dimension this whole exercise is pointless; the traditional Flow-box
Theorem applies to merely continuous vector fields. If $f:\mathbb{R\rightarrow
R}$ is continuous and $f\left(  x_{1}\right)  \neq0$ then there exists a
diffeomorphism $\phi$ which transfers $f$ near $x_{1}$ to the constant vector
field $g=1$. Explicitly,%
\[
\phi\left(  x\right)  =\int_{x_{1}}^{x}\frac{1}{f\left(  t\right)  }dt
\]
which is a diffeomorphism near $x_{1}$ since $f\left(  x_{1}\right)  \neq0$.
Therefore%
\begin{align*}
\phi^{\prime}\left(  x\right)  \cdot f\left(  x\right)   &  =1\\
\text{so\hspace{0.3in}}\left(  d\phi\right)  _{x}\left(  f\left(  x\right)
\right)   &  =g\left(  \phi\left(  x\right)  \right)
\end{align*}
and so $\phi$ transfers $f$ to $g$.
\end{remark}

\begin{example}
Curiously, there exist discontinuous vector fields which can be transferred to
constant vector fields by a lipeomorphism. Consider $f,g:\mathbb{R\rightarrow
R}$ where
\[
f\left(  x\right)  =\left\{
\begin{array}
[c]{cc}%
1 & x<1\\
2 & x\geq1
\end{array}
\right.
\]
and $g\left(  x\right)  =1$. Then $f$ is locally topologically conjugate to
$g$ via
\[
\phi\left(  x\right)  =\left\{
\begin{array}
[c]{cc}%
x & x<1\\
\frac{x+1}{2} & x\geq1\text{.}%
\end{array}
\right.
\]
\end{example}

\end{document}